\documentclass[11pt]{amsart}
\usepackage{graphicx} 
\usepackage{lineno}

\usepackage{amsthm}
\usepackage{amssymb}
\usepackage{amsfonts}
\usepackage{mathtools}
\usepackage{mathrsfs}
\usepackage{graphicx}
\usepackage{color}
\usepackage{enumitem}
\usepackage{setspace}
\usepackage[english]{babel}
\usepackage{csquotes}
\usepackage{etoolbox}
\usepackage{thmtools}
\usepackage{yfonts}
\usepackage{mfirstuc}
\usepackage{textgreek} 

\usepackage{amsmath}
\usepackage{amsthm, amssymb, amsfonts}
\usepackage{graphicx, color}
\usepackage{caption, subcaption} 
\usepackage{multicol} 
\usepackage{stmaryrd}
\usepackage{tikz-cd} 
\usepackage{lipsum}  
\usepackage[capitalise,noabbrev, nameinlink]{cleveref} 
\Crefname{part}{Part}{Parts}
\Crefname{step}{Step}{Steps}
\Crefname{prop}{Proposition}{Propositions}
\Crefname{prob}{Problem}{Problems}
\newcommand{\forces}[1]{\Vdash\text{\say{#1}}}

\usepackage{enumitem}
\renewcommand{\leq}{\leqslant}
\renewcommand{\geq}{\geqslant}

\renewcommand{\epsilon}{\varepsilon}

\newcommand{\dom}{\text{dom}}
\newcommand{\stickT}{%
\setbox255=\hbox{\raise1ex\hbox{$\hspace{0.2pt}\,\bullet\,$}}
\mathord{\rlap{\hbox to\wd255{\hss\hbox{$|$}\hss}}
\box255}
}
\newcommand{\stickS}{%
\setbox255=\hbox{\raise0.6ex\hbox{$\scriptstyle\bullet$}}
\mathord{\rlap{\hbox to\wd255{\hss\hbox{$\scriptstyle|$}\hss}}
\box255}
}

\numberwithin{equation}{section}
\theoremstyle{plain}
\newtheorem{theorem}[equation]{Theorem}

\newtheorem{proposition}[equation]{Proposition}
\newtheorem{corollary}[equation]{Corollary}

\theoremstyle{definition}
\newtheorem{definition}[equation]{Definition}

\newtheorem{problem}[equation]{Problem}

\theoremstyle{remark}

\usepackage{mathtools}
\usepackage{amsfonts}
\usepackage{amssymb}
\usepackage{amsthm}
\usepackage{amsmath}
\usepackage{dirtytalk}
\usepackage{mathrsfs}  
\usepackage{lineno}

 \usepackage{multicol}  
 \usepackage{bbm}
 \usepackage{tabularx}
\newenvironment{claimproof}[1][\proofname]
{\begin{proof}[#1]}
{\end{proof}}

\title{$\kappa$-Barely independent families and Tukey types of ultafilters}
\author{Jorge Cruz}

\begin{document}
\maketitle
\begin{abstract} Given two infinite cardinals $\kappa$ and $\lambda$, we introduce and study the notion of a $\kappa$-barely independent family over $\lambda.$ We provide some conditions under which these types of families exist. In particular, we relate the existence of large $\kappa$-barely independent families with the generalized reaping numbers $\mathfrak{r}(\kappa,\lambda)$ and use these relations to give conditions under which every uniform ultrafilter over a given cardinal $\lambda$ is both Tukey top and has maximal character. Finally, we show that $\mathfrak{p}>\omega_1$  the non-existence of  barely independent families over $\omega_1.$
\end{abstract}

\section{Introduction}
\footnote{Mathematics subject classification: 03E17, 03E10, 03E35, 03E05, 03E20}
Given a set $X$, an \emph{independent family} over $X$ is a family $\mathcal{I}$ of subsets of $X$ so that $|\bigcap \mathcal{A}\backslash \bigcup \mathcal{B}|=|X|$ whenever $\mathcal{A},\mathcal{B}\subseteq \mathcal{I}$ are disjoint disjoint and finite. Independent families first appeared in the work of Hausdorff. Since then, they have become central objects in Set theory and Topology. For example, regarding the theory of ultrafilters,  Hajnal and Juh\'asz used the existence of an independent family of size $2^\kappa$ over $\kappa$ (here, $\kappa$ is an infinite cardinal) to prove that there are $2^{2^\kappa}$ many ultrafilters over $\kappa$ of character $2^\kappa$. Similarly in \cite{CategorycofinaltypesII}, Isbell used large independent families in order to construct \emph{Tukey Top} ultrafilters over $\kappa$ (defined below). Given a set $X$ and  $A\subseteq X$, let $A^0=A$ and $A^1=X\backslash A$. The results stated above come from the fact that if  $\mathcal{I}$ is an independent family over $\kappa$ and $f:\mathcal{I}\longrightarrow 2$ is an arbitrary function, there is at least one ultrafilter $\mathcal{U}$ over $\kappa$ extending both the family $$\mathcal{I}^f:=\{A^{f(A)}\,:\,A\in \mathcal{I}\}$$
    and the family
    $$\{ \kappa\backslash (\bigcap\limits_{A\in \mathcal{A}} A^{f(A)})\,:\,\mathcal{A}\in [\mathcal{I}]^{\omega}\}.$$ 

It is easy to see that the character of any such ultrafilter is at least $|\mathcal{I}|$ and furthermore, it is Tukey above $[ \,|\mathcal{I}|\,]^{<\omega}$ aIn particular, if the size of $\mathcal{I}$ is $2^\kappa$, then such ultrafilter is both Tukey top an has maximal character. It is easy to see that the converse is also true, that is, if an ultrafilter over $\kappa$ is Tukey above a cardinal $\lambda$, then there is an independent family $\mathcal{I}$ witnessing such fact (in the sense described above). 

One of the questions guiding this research is under which conditions,for a given $\kappa$ and $\lambda$, there exists single independent family $\mathcal{I}$  of size $\lambda$ witnessing that every uniform ultrafilter over $\kappa$ is Tukey above $[\lambda]^{<\omega}$. 
In the process of answering this question, we arrive to the concept of \emph{barely independent families} and the more general notion of \emph{$\kappa$-barely independent families}. It turns out that this kind of families may only exist on uncountable cardinals (see Proposition \ref{nobarelyomega}). Hence, this topic also fits into the increasing line of research of \say{higher independent families} in the uncountable setting such as \emph{$\sigma$-independent families}, \emph{$(\theta,\tau)$-independent families}, \emph{strongly independent families}, $\mathcal{C}$-independent families, etc., which can be found in papers such as \cite{maximalsigmaindependent}, \cite{independentfamiliesresolvability}, \cite{higherindependence}, \cite{generalizedindependence} and \cite{strongindependencespectrum}. The main difference between those concepts and this new one, is that while they all require  large boolean combinations to be always large in some sense, we require large boolean combinations to be always small. 

As we mentioned, the concept we are defining is in its core related to the Tukey order and the character of ultrafilters. For that reason, it is convenient to give a quick review of these topics before we can state the main results of the paper.

Given an ultrafilter $\mathcal{U}$, a \emph{base} for $\mathcal{U}$ is a family $\mathcal{B}\subseteq \mathcal{U}$ which which generates $\mathcal{U}$. That is, every element of $\mathcal{U}$ contains an element of $\mathcal{B}$. The character of an ultrafilter, written as $\chi(\mathcal{U})$, is the minimal cardinality of one of its bases. If $\mathcal{U}$ is an ultafilter over some set $X$, we say that it is \emph{uniform}, provided that each of its elements has size $|X|$. It is easy to see that each uniform ultrafilter over a regular cardinal $\kappa$ has character at least $\kappa$ and at most $2^{\kappa^+}$. For the case of $\omega$ it is known that the existence of an ultrafilter of character strictly smaller than $\mathfrak{c}$ is independent from $ZFC$. Suprisingly, for $\omega_1 $ this question is still \say{very} open.
\begin{problem}[Kunen] Is it consistent that there is a ultrafilter over $\omega_1$ whose character is smaller than $2^{\omega_1}$?
\end{problem}
Let $D$ and $E$ be two directed sets. $D$ is said to be \emph{Tukey below} $E$, written as $D\leq_T E$, provided that there is a function $f:D\longrightarrow E$ mapping unbounded sets of $D$ to unbounded sets of $E$. If both $D\leq_T E$ and $E\leq_T D$ occur, we say that $D$ is \emph{Tukey equivalent} to $E$ and write is as $D\equiv_T E$. As proved by Tukey in \cite{Tukey} when he first introduced this concept, $D\equiv_T E$ is equivalent to the existence of a directed set which embeds cofinally in both $D$ and $E$.

For more than 80 years, Tukey theory has been an active area of research, particularly among set theorists and topologists.  A lot of effort has been put into developing classification theories for classes of structures in terms of the Tukey order. For example, in \cite{Tukey} it was already shown that the partial orders $1, \omega,\omega_1,\omega\times \omega_1$ and $[\omega_1]^{<\omega}$ are non Tukey equivalent when equipped with their natural orderings. In the celebrated \cite{directedsetscofinaltypes}, Todor\v{c}evi\'{c} proved that under PFA, any other partial order of size $\aleph_1$ is Tukey equivalent to one of the five previously mentioned orders. This result is known to be false under several distinct assumptions (see \cite{CategorycofinaltypesII}, \cite{sevencofinaltypes} and \cite{PartitionProblems}). We have the so called \say{Isbell problem}. In \cite{CategorycofinaltypesII}, Isbell constructed a \emph{Tukey top} ultrafilters over each regular cardinal $\kappa$, that is an ultrafilter which is maximal among all directed sets of size $2^\kappa$ when ordered by reverse inclusion. In that same paper, he asked:
\begin{problem}[Isbell, Problem 2, \cite{CategorycofinaltypesII}] How many non-equivalent ultrafilters  exist on a set of power $\aleph_0$?  
\end{problem}
Several contributions have been made to this question from which we highlight \cite{tukeytypesofultrafilters} and \cite{surveytukeytheory}. Recently in \cite{ontheisbellproblem}, Cancino and Zapletal finally settled Isbell problem by constructing a model of ZFC in which every non-principal ultrafilter over $\omega$ is Tukey top. Now, the next natural step is to settle what we call here \say{the generalized Isbell problem}:
\begin{problem}Let $\kappa$ be a cardinal. How many non-equivalent ultrafilters exist on a set of power $\kappa$?
\end{problem}
In \cite{cofinaltypesmeasurable}, Benhamou and Dobrinen considered this question for the case where $\kappa$ is a measurable cardinal. Recently, Benhamou, Moore and Serafin  \cite{thecofinalityofultrafiltersomega1} have settled the question for $\kappa=\omega_1$ by showing its independence from ZFC.\\

In this paper we will focus in analyzing the existence and non-existence of ($\kappa$-)barely independent families under several assumptions. The main contributions are the following:
\begin{itemize}
    \item In Theorem \ref{independentsuccesor} we show that if $2^\kappa=\kappa^+$, then there is a $\kappa^+$-barely independent family over $\kappa^+$ of size $\kappa^+$. On the other hand, there are no $\omega$-barely independent families over $\omega$ (see Proposition \ref{nobarelyomega}). These two results provide a good contrast between the concept of $\kappa$-barely independence in countable and uncountable ordinals.
    \item We show that there is a strong relation between the cardinals $\mathfrak{r}(\kappa,\lambda)$ as defined in \cite{invariantsuncountablecardinals} and $\kappa$-barely independent families. Namely, in Theorem \ref{rimpliesbarely} we show that if there is a $\kappa$-barely independent family over $\mu$ of size $\lambda$, then $\mathfrak{r}(\kappa,\lambda)\geq cof(\mu)$. Conversely, in Theorem \ref{barelyimpliesr} we show that if $\mathfrak{r}(\kappa,\lambda)=2^\kappa$ and $\lambda^{<\kappa}\leq cof(2^\kappa)$, then there is a $\kappa$ barely independent family over $2^\kappa$ (also over $cof(2^\kappa)$) of size $\lambda.$
    \item In Theorem \ref{teoprinma}, we show that there are no barely independent families over any uncountable cardinal smaller than $\mathfrak{p}$. 
    \item Lastly, we will use the concept of barely independent families, not only to show that for any regular cardinal $\kappa$ there is a cardinal preserving forcing extension of the universe in which every uniform ultrafilter over $\kappa$ is Tukey top (something which is also proved in\cite{thecofinalityofultrafiltersomega1}) and has maximal character, but for uncountable cardinals between $\omega$  and $\mathfrak{c}$ we will be able to axiomatize this result in terms of cardinal arithmetic and cardinal invariants of the continuum. Specifically, we provide the following contribution to both the Kunen problem and the generalized Isbell problem:\\
\begin{center}If $2^{cof(\mathfrak{c})}=\mathfrak{c}=\mathfrak{r}(\omega,cof(\mathfrak{c}))$, then every uniform ultrafilter over $cof(\mathfrak{c})$ is Tukey top and has character $2^\kappa$.
\end{center}
It is important to remark that this equality holds in any generic extension obtained after adding $\mathfrak{c}$-many Cohen reals. 
\end{itemize}
\vspace{0.3cm}

The structure of the paper is as follows: In \emph{Section} \ref{notationsection} we briefly recall some important notation. In \emph{Section} \ref{barelyindependentsection} we introduce $\kappa$-barely independent fammilies and prove some first results. In \emph{Section} \ref{cardinalcharacteristicssection} we relate $\kappa$-barely independent families with the cardinals $\mathfrak{r}(\kappa,\lambda)$ and show that barely independent families over $\kappa$ exist after adding $\kappa$ many Cohen reals. As a corollary of this, we will also be able to show that consistently there are  barely independent families over $\omega_1$, but not over $\mathfrak{c}$. Finally, in Section \ref{openproblemssection}, we pose some open problems.\\

The author would like to thank Cancino, Guzm\'an, Hru\v{s}ak and Todor\v{c}evi\'c for their helpful comments during the preparation of this work.
\section{Notation}\label{notationsection}
The notation and terminology used here mostly standard and it follows \cite{schemescruz}.  Given a set $X$ and a (possibly finite) cardinal $\kappa$, $[X]^\kappa$ denotes the family of all subsets of $X$ of cardinality $\kappa$. The sets $[X]^{<\kappa}$ and $[X]^{\leq \kappa}$ have the expected meanings.  $\mathscr{P}(X)$ denotes the power set of $X$. For a set of ordinals $X$, we denote by $ot(X)$ its order type. We identify $X$ with the unique strictly increasing function $h:ot(X)\longrightarrow X$. In this way, $X(\alpha)=h(\alpha)$ denotes the $\alpha$ element of $X$ with respect to its increasing enumeration. Analogously, $X[A]=\{X(\alpha)\,:\,\alpha\in A\}$ for $A\subseteq ot(X)$. For a function $h$, its domain its denoted as $\dom(h)$. If $h$ is a partial function from a set $X$ to a set $Y$, we denote it as $h;X\longrightarrow Y$. The set $Fin(X,2)$ denotes the set of all finite partial functions from $X$ to $2$. This set ordered with the reverse inclusion is equivalent to the Cohen forcing for adding $|X|$-many Cohen reals. In the particular case in which $X=Y\times Z$ and $p\in Fin(X,2)$, we define $dom_0(p):=\{y\in Y\,:\,\exists z\in Z\,((y,z)\in \dom(p))\}$ . The set $dom_1(p)$ is defined in an analogous way. As mentioned before, if $A\subseteq X$, then $A$ is denoted as $A^0$ and $X\backslash A$ as $A^1$. Given $\mathcal{A}\subseteq \mathscr{P}(X)$ and $h:\mathcal{A}\longrightarrow 2$, $\mathcal{A}^h$ stands for the set $\{A^{h(A)}\,:\,A\in \mathcal{A}\}.$ A family $\mathcal{D}$ is called a \emph{$\Delta$-system} with root $R$ if $|\mathcal{D}|\geq 2$ and $X\cap Y= R$ whenever $X,Y\in \mathcal{D}$ are different. 

Given a forcing notion $\mathbb{P}$, we denote by $\mathfrak{m}(\mathbb{P})$ its Martin's number, i.e. , the minimal cardinal $\kappa$ for which there are $\kappa$-many dense sets over $\mathbb{P}$ which are not intersected by any filter. The Martin's number $\mathfrak{m}$ is defined as the minimum of all the Martin's numbers of ccc forcings. For us, MA is the statement \say{$\mathfrak{m}=\mathfrak{c}>\omega_1$}.

\section{Barely independent families}\label{barelyindependentsection}
Assume that $\lambda$ is an infinite cardinal. When is it true that for a given independent family $\mathcal{I}$, every uniform ultrafilter over $\lambda$ contains both $\mathcal{I}^h$ and $co^h_\omega(\mathcal{I}):=\{\lambda\backslash (\bigcap\limits \mathcal{A})\,: \mathcal{A}\in [\mathcal{I}^h]^\omega\}\,$?
Evidently, every ultrafilter extends $\mathcal{I}^{h_\mathcal{U}}$ for a unique $h_\mathcal{U}$. So the problem lies in the second family. If the intersection of $\mathcal{A}$ has cardinality less than $\kappa$ for each countable $\mathcal{A}\subseteq \mathcal{I}^h$,the uniformity of $\mathcal{U}$ would imply that it extends $co^h_\omega(\mathcal{I})$. On the other hand, if for some $\mathcal{A}\in [\mathcal{I}]^h$ and $h':\mathcal{A}\longrightarrow 2$ it happens that $\bigcap \mathcal{A}^{h'}$ has size $\lambda$ and $\mathcal{U}$ is a normal ultrafilter having such intersection, then $h_\mathcal{U}$ will extend $h'$ and $co^{h_\mathcal{U}}_\omega(\mathcal{I})$ will not be contained in $\mathcal{U}$.  Thus, independent families having the property discussed above are exactly the ones for which all of  their infinite boolean combinations have size less than $\lambda$. In this way, we arrive to the notion of barely independent families.

\begin{definition}[Barely independent families] Let $\kappa\leq \lambda$ be a cardinal and $\mathcal{I}\subseteq [\lambda]^\lambda$. We say that $\mathcal{I}$ is \emph{$\kappa$-barely independent} provided that $|\mathcal{I}|\geq \kappa$ and:\begin{itemize}
    \item $|\bigcap \mathcal{A}^h|=\lambda$ for each $\mathcal{A}\in [\mathcal{I}]^{<\kappa}$ and $h:\mathcal{A}\longrightarrow 2$.
    \item $|\bigcap \mathcal{A}^h|<\lambda$ for each $\mathcal{A}\in [\mathcal{I}]^\kappa$ and $h:\mathcal{A}\longrightarrow 2.$
    \end{itemize}
In the particular case in which $\kappa=\omega$, we simply call $\mathcal{I}$ \emph{barely independent}.
\end{definition}
As we remarked at the begining of this paper, part of the importance of properties which we are studying lie in their relation with the Tukey type and the character of ultrafilters. We capture those ideas in the next two propositions.
\begin{proposition}Let $\kappa\leq \lambda$ be two cardinals and assume that there is a $\kappa$-barely independent family $\mathcal{I}$ in $\lambda$. Then every uniform ultrafilter in $\lambda$ is Tukey above $[\,|\mathcal{I}|\,]^{<\kappa}$.
    
\end{proposition}

\begin{proposition}Let  $\lambda$ be a cardinal and assume that there is a barely independent family $\mathcal{I}$ in $\lambda$. Then every uniform ultrafilter in $\lambda$ has  character greater or equal than $|\mathcal{I}|$.   
\end{proposition}
The concept of barely independence only makes sense in the uncountable setting. In particular, our technics do not yield any new information regarding the Isbell problem over $\omega$. Except maybe showing that one of most naive ways of trying to solve such problem is far from working.

\begin{proposition}\label{nobarelyomega}There are no barely independent families in $\omega$.
\begin{proof} Let us assume towards a contradiction that $\mathcal{I}\subseteq [\omega]^{<\omega}$ is a barely independent family. \\\\
\noindent
\underline{Claim 1:} $\{\min(A)\,:\,A\in \mathcal{I}\}$ is bounded.
\begin{claimproof}[Proof of claim] Otherwise, we can recursively build an infinite sequence $\langle A_n\rangle_{n\in \omega}$ of  elements of $\mathcal{I}$ so that $\min(A_{n+1})> \min((\bigcap\limits_{i\leq n}A^1_i)\backslash \min (A_n))$  for each $n$. In particular,  the infinite set $$\{\min\big((\bigcap\limits_{i\leq n}A^1_i)\backslash \min (A_n)\big)\,:\,n\in \omega\}$$ is contained in $\bigcap \limits_{i\in  \omega}A^1_i$. This contradicts the fact that $\mathcal{I}$ is  barely independent. 
\end{claimproof}

We will now arrive to a contradiction by showing that there is an infinite $\mathcal{A}\in \mathcal{I}$ for which $\bigcap \mathcal{A}=\bigcap\limits_{A\in \mathcal{A}}A^0$ is infinite. According to Claim 1 above, there is $n_0\in \omega$  so that $\mathcal{I}_0=\{A\in \mathcal{I}\,:\,\min(A)=n_0\}$ is  infinite. Let $A_0\in \mathcal{I}_0$. It is straightforward that the  family $\mathcal{I}'_0=\{(A_0\backslash (n_0+1))\cap A\,:\,A\in \mathcal{I_0}\backslash \{A_0\}\}$ is a barely independent family in $A_0\backslash (n_0+1)$. By applying Claim 1 to this family, we can find $n_1\in A_0\backslash n_0$ for which $\mathcal{I}_1=\{A\in \mathcal{I_0}\,:\,\min (A\cap (A_0\backslash (n_0+1)))=n_1\}\,$ is infinite. By continuing this process we can construct an increasing sequence $\langle n_i\rangle_{i\in \omega}$ of natural numbers and a sequence $\langle A_i\rangle_{i\in \omega}$ of elements of $\mathcal{I}$ in such way that $$\min\big(A_k\cap ((\bigcap_{i\leq j} A_i\backslash ( n_j+1)\big)=n_{j+1}$$ for any two $j<k$. In particular, we have that $\{n_i\,:\,i\in \omega\}\subseteq \bigcap\limits_{i\in \omega}A_i$ which is the desired contradiction.
\end{proof}  
\end{proposition}

As it is shown below, the study of  $\kappa$-barely independent families is only possible for cardinals whose cofinality is at most than $2^\kappa.$

\begin{proposition}Let $\kappa$ be a cardinal and assume that there is a $\kappa$-barely independent family $\mathcal{I}$ over $cof(\lambda)$. Then there is a barely independent family over $\lambda$ of the same size.
\begin{proof}Indeed, let $h:\lambda\longrightarrow cof(\lambda)$ be such that $\langle|h^{-1}[\{\alpha\}]\rangle_{\alpha\in cof(\lambda)}$ converges to $\lambda.$ It should be clear that $\mathcal{I}'=\{h^{-1}[A]\,:\,A\in \mathcal{A}\}$ works. 
\end{proof} 
\end{proposition}

\begin{proposition}Let $\kappa,\lambda$ be two regular cardinals so that $cof(\lambda)>2^\kappa$. There are no barely $\kappa$-independent families over $\lambda$.
\begin{proof}Let $\mathcal{I}$ be an independent family over $\lambda$.  We will show that $\mathcal{I}$  is not $\kappa$-barely independent. Without loss of generality we may assume that $\mathcal{I}$  has size $\kappa.$ From this fact, it is easily seen that there is a partition $\{  B_\alpha\,:\,\alpha\in \lambda
\}$ of $\lambda$ into sets of size $\kappa$ in such way that the family $\mathcal{I}_\alpha=\{A\cap B_\alpha\,:\,A\in \mathcal{I}\}$ is an independent family over $B_\alpha$ for each $\alpha\in \lambda$. Given any such $\alpha,$ let $\varphi_\alpha:B_\alpha\longrightarrow \kappa$ be bijection and define $f_\alpha:\mathcal{I}\longrightarrow [\kappa]^{\kappa}$ as: $$f_\alpha(A)=\varphi_\alpha(A\cap B_\alpha).$$
Since $\mathcal{I}$ is  has size $\kappa$ and $cof(\lambda)>2^\kappa$ , we can find $X\in [\lambda]^{\lambda}$ and $f:\mathcal{I}\longrightarrow \kappa$ such that $f_\alpha=f$ for each $\alpha\in X$. It is straightforward that $\{f[A]\,:\,A\in \mathcal{I}\}$ is independent. In particular there is $h:\mathcal{I}\longrightarrow 2$ so that $\bigcap\limits_{A\in \mathcal{I}} f(A)^{h(A)}$ is non-empty. From this it follows that $(\bigcap \mathcal{A}^h)\cap B_\alpha\not=\emptyset$ for any $\alpha \in X$. Since $X$ has size $\lambda$, so does $\bigcap \mathcal{I}^h$.
\end{proof}
\end{proposition}

Even though there are no $\omega$-barely independent families over $\omega$, there are (consistently) $\mu$-barely independent families over $\mu$ whenever $\mu$ is a succesor cardinal.

\begin{theorem}\label{independentsuccesor}Let $\kappa$ be a cardinal such that $2^{k}=k^+$. Then there is an $\kappa^+$-barely independent family over $k^+$.
\begin{proof}Let $\mathcal{F}$ be the set of all partial functions from $k^+$ into $2$ whose domain has size at most $\kappa$ and let $\hat{h}:\kappa^+\longrightarrow \mathcal{F}$ be a surjective function with cofinal fibers. Let us consider a club $C$ in $\kappa^+$ consisting of ordinals greater than $\kappa$ in such way that for $\alpha,\beta\in \kappa^+$, if $\alpha< C(\beta)$, then $\dom(\hat{h}(\alpha))\subseteq C(\beta).$ Now, we enumerate $[\kappa^+]^{\kappa}$ as $\langle R_\alpha\rangle_{\alpha\in \omega_1}$ in such way that $R_\alpha\subseteq C(\alpha)$ for each $\alpha$. Note that since $\{R_\xi\,:\,\xi\leq \alpha \}$ is a $\kappa$-sized family of  subsets of $\alpha$, each of cardinality $\kappa$ , we can fix a set $S_\alpha\subseteq C(\alpha)$ which splits $R_\xi$ for each $\xi\leq\alpha$. Finally, let $\mu:\kappa^+\longrightarrow \kappa^+$ be defined as $$\mu(\beta)=\max( \gamma\in \kappa^+\,:\,C(\gamma)\leq \beta).$$
Note that $\mu$ is well-defined since $C$ is a club. Given $\beta\in \kappa^+$, we put $$H_\beta=\{\alpha \in \kappa^+\,:\, \beta\in \dom(\hat{h}(\alpha))\text{ and }\hat{h}(\alpha)(\beta)=0\}$$
and define $A_\beta:=S_{\mu(\beta)}\cup H_\beta$. 
In the following two claims, we will prove that the family $\mathcal{I}:=\{A_\beta\,:\,\beta\in \kappa^+\}$ is $\kappa^+$-barely independent.\\\\
\noindent
\underline{Claim 1}: Let $\mathcal{A}\in [\omega_1]^{<\omega_1}$ and $h:\mathcal{A}\longrightarrow 2.$ Then $\bigcap\limits_{\beta\in \mathcal{A}}A^{h(\beta)}_\beta$ has size $\kappa^+$.
\begin{claimproof}[Proof of claim] Indeed, If $\alpha\in \omega_1$ is such that $\hat{h}(\alpha)=h$, then $$\alpha\in \bigcap\limits_{\beta\in \mathcal{A}} H^{h(\beta)}_\beta\subseteq \bigcap\limits_{\beta\in \mathcal{A}} A_\beta^{h(\beta)}.$$
Since each fiber of $\hat{h}$ is unbounded in $\kappa^+$, this finishes the proof.
\end{claimproof}
\noindent
\underline{Claim 2}: Let $\mathcal{A}\in [\omega_1]^{\omega_1}$ and $h:\mathcal{A}\longrightarrow 2$. Then $\bigcap\limits_{\beta\in \mathcal{A}} A_\beta^{h(\beta)}$ has size strictly less than $\kappa$.
\begin{claimproof}[Proof of claim] Assume towards a contradiction that such intersection has size at least $\kappa$ and let $\xi\in \kappa^+$ be such that $$R_\xi\subseteq \bigcap\limits_{\beta\in \mathcal{A}} A_\beta^{h(\beta)}. $$Consider  $\beta\in \mathcal{A}$ for which  $\mu(\beta)>\xi$. Then $S_{\mu(\beta)}$ splits $R_\xi$. Furthermore, $A_\beta\cap C(\mu(\beta))=S_{\mu(\beta)}$. Hence neither $A_\beta^0$ nor $A_\beta^0$ contain $R_\xi$. This contradicts the fact that $R_\xi \subseteq A_\beta^{h(\beta)}$, so we are done.
\end{claimproof}
    
\end{proof}   
\end{theorem}
In \cite{higherindependence}, a \emph{strongly independent family} $\mathcal{I}$ (over a cardinal $\kappa)$ is defined as an independent family over $\kappa $ such that for any $\mathcal{A}\in [\kappa]^{<\kappa}$ and each $h:\mathcal{A}\longrightarrow 2$, the set $\bigcap \mathcal{A}^h$ has size $\kappa$. Evidently every $\kappa$-barely independent family over $\kappa$ is strongly independent. Hence, Theorem \ref{independentsuccesor} slightly generalizes Theorem 1.10 from \cite{generalizedindependence}.
\begin{corollary}\label{corostronglyindependent}Let $\kappa$ be an infinite cardinal. The followings are equivalent:
\begin{enumerate}
    \item $2^{\kappa}=\kappa^+$.
    \item  There is a $\kappa^+$-barely independent family over $\kappa^+.$
    \item There is a strongly independent family of size $\kappa^+$ over $\kappa^+.$
    \item There is a strongly independent family of size $\kappa$ over $\kappa^+.$
\end{enumerate}  
\end{corollary}

\section{Barely independent families and cardinal characteristics of the continuum}\label{cardinalcharacteristicssection}

Let $R$ and $S$ and be two infinite sets. We say that $S$ \emph{splits} $R$  in case both $R\cap S$ and $R\backslash S$ are infinite. If this do not occur, we say that $R$ \emph{reaps} $S$. That is, either $R\subseteq^* S$ or $R\subseteq ^*X\backslash S$. In \cite{invariantsuncountablecardinals}, the number $\mathfrak{r}(\kappa,\lambda)$ is defined as the minimum size of a family $\mathcal{R}\subseteq [\lambda]^{\kappa}$ so that every infinite $S\subseteq \kappa$ is reaped by at least one member of $\mathcal{R}$. We call any such family a \emph{reaping family}. Of course, here $\kappa$ and $\lambda$ stand for two distinct cardinals. It should be clear that $\mathfrak{r}(\kappa,\kappa)$ is equal to the cardinal $\mathfrak{r}_\kappa$ studied in papers such as and \cite{cardinalcharacteristicskappasmallu}, \cite{remarksgeneralizedUDR}, \cite{Galvinghajnalgeneralizedcardinal} and \cite{tworesultsinvariantsuncountable}. In particular,  $\mathfrak{r}(\omega,\omega)$ is the same as the classical reaping number $\mathfrak{r}$. The next proposition is well-known and easy to prove.

\begin{proposition}let $\kappa_0,\kappa_1\leq \lambda_0,\lambda_1$ be infinite cardinals. The following properties hold:
\begin{itemize}
    \item If $\kappa_0\leq \kappa_1$, then $\mathfrak{r}(\kappa_0,\lambda_0)\leq \mathfrak{r}({\kappa_1,\lambda_0}).$
    \item If $\lambda_0\leq \lambda_1$, then $\mathfrak{r}(\kappa_0,\lambda_0)\geq \mathfrak{r}(\kappa_0,\lambda_1).$ 
    \item $\kappa_0^+\leq \mathfrak{r}(\kappa_0,\lambda_0)\leq \lambda^{\kappa_0}\leq 2^{\lambda_0}.$
    \end{itemize}
Particularly, this means that for a fixed $\kappa$, we have the following chain of inequalities: $$\kappa^+\leq \mathfrak{r}(\kappa,2^\kappa)\leq \dots\leq \mathfrak{r}(\kappa,\kappa^+)\leq \mathfrak{r}(\kappa,\kappa)\leq 2^\kappa.$$
\end{proposition}
In \cite{higherindependence}, Fischer and Montoya observed that, just as in the countable setting,  if a strongly independent family $\mathcal{I}$ over a regular cardinal $\kappa$ is such that $|\mathcal{I}|\leq \mathfrak{r}(\kappa,\kappa)$, then $\mathcal{I}$ is not maximal. We will now show that there is also a relation between reaping numbers and the existence of $\kappa$-barely independent families.

\begin{theorem} \label{barelyimpliesr}
If there is a $\kappa$-barely independent family over $\mu$ of size $\lambda$, then $\mathfrak{r}(\kappa,\lambda)\geq cof(\mu)$. 
\begin{proof}
Let $\mathcal{R}\subseteq [\lambda]^{\kappa}$ by a family of size less than $cof(\mu)$. We will prove that $\mathcal{R}$ is not reaping. For this, we fix a $\kappa$-barely independent  family $\mathcal{I}=\{A_\alpha\,:\,\alpha\in \lambda\}$  over $\mu.$   Given $R\in \mathcal{R}$, $i\in 2$ and $a\in[R]^{<\omega}$, let $h_{R,a,i}:R\backslash a\longrightarrow 2$ be the constant function with value $i$. Since $\mathcal{I}$ is $\kappa$-barely independent, there is $\alpha_{R,a,i}<\mu$ such that $$\bigcap \limits_{\xi\in R\backslash a}A_\xi^{h_{R,a,i}(A)}\subseteq \alpha_{R,a,i}.$$ As $|\mathcal{R}|<cof(\mu)$, there exists some $\alpha$ so that $\alpha>\alpha_{R,n,i}$ for all  $R\in \mathcal{R}$, $a\in[R]^{<\kappa}$ and $i\in 2$. Let $S=\{\xi\in \lambda\,:\,\alpha\in A_\xi\}$ and fix $R\in \mathcal{R}$. We claim that $S$ splits $R$. Indeed, given $a\in [R]^{<\kappa}$ we have that $$\bigcap \limits_{\xi\in R\backslash a}A_\xi^1 \subseteq \alpha_{R,n,1}\subseteq \alpha.$$ In particular, there is $\xi\in R\backslash a$ for which $\alpha\not\in A_\xi^1=\kappa\backslash A_m$. By definition, $\xi\in R\cap S$. This proves that $R\cap S$ is infinite. Analogously, using the fact that $\alpha_{R,n,0}<\alpha$ we can prove that $R\backslash S$ is infinite as a well. Since $R$ was arbitrary, we conclude that $\mathcal{R}$ is not reaping.
\end{proof}
\end{theorem}
\begin{theorem}\label{rimpliesbarely} Let $\kappa$ be a regular cardinal,  $\mu=cof(2^\kappa)$ and $\lambda\leq \mu$ be such that $\lambda^{<\kappa}\leq \mu$. If $\mathfrak{r}(\kappa,\lambda)=2^\kappa$, then there is a $\kappa$-barely independent family over $\mu$ of size $\lambda$.
\begin{proof}Since $|[\lambda]^\kappa|=2^\kappa$ and $\kappa$ is regular, we can write $[\lambda]^\kappa$ as an increasing union of a sequence $\langle \mathcal{R}_\alpha\rangle_{\alpha\in \mu}$ satisfying the following properties for each $\alpha\in \mu$:
\begin{itemize}
\item $|\mathcal{R}_\alpha|<2^\kappa$.
\item For any $R\in \mathcal{R}_\alpha$ and each $X\in [R]^{<\kappa}$, there is $R'\in[R]^\kappa$ such that $X\cap R=\emptyset.$
\end{itemize}

Since $\lambda^{<\kappa}\leq \mu$, we can also fix $\sigma:\mu\longrightarrow \{h;\lambda\longrightarrow 2\,:\,|\dom(h)|<\kappa\}$  surjective with cofinal fibers. Given $\alpha\in \mu$, the equality $r(\kappa,\lambda)=2^\kappa$ implies that $\mathcal{R}_\alpha$ is not  reaping. In other words, there is $S_\alpha\subseteq\lambda$ which splits every element of $\mathcal{R}_\alpha$. Let $\alpha\in \lambda$. We define $A_\alpha\subseteq \mu$ as follows: $$A_\alpha:=\{\xi\in \mu\,:\,(\alpha\in \dom(\sigma(\xi))\textit{ and }\sigma(\xi)(\alpha)=0)\text{ or }(\alpha\notin \dom(\sigma(\xi))\text{ and }\alpha\in S_\xi)\}.$$
We claim that $\mathcal{I}:=\{A_\alpha\,:\,\alpha\in \lambda\}$ is $\kappa$-barely independent. Indeed, for any  partial function $h;\lambda\longrightarrow 2$ such that $|\dom(h)|<\kappa$, it happens that $$\sigma^{-1}[\{h\}]\subseteq \bigcap\limits_{\alpha\in \dom(h)}A_\alpha^{h(\alpha)}.$$
Since the function $\sigma$ has cofinal fibers, we conclude that the intersection written above has cardinality $\mu$. In order to finish,  take a  $B\subseteq \lambda$  of size $\kappa$ and $h:B\longrightarrow 2$. Note that  there is $i\in 2$ so that $R:=h^{-1}[\{i\}]$ has size $\kappa$. For such $i$ we know that there is $\alpha \in \mu$ such that $R\in \mathcal{R}_\alpha$. By construction, given $\beta\geq \alpha$ there is $R'\in \mathcal{R}_\beta$ which is contained in $R$ and disjoint from $\dom(\sigma(\beta))$. Note that  $S_\beta$ splits $R'$. Thus, there is $\xi\not\in \dom(\sigma(\beta))$ so that $\xi\in R\cap S_\beta^{1-i}$. This means that $\beta\not\in A_\xi^i=A_\xi^{h(\xi)}$. In particular we have that $$\beta\not\in \bigcap \limits_{\zeta\in \dom(h)}A_\zeta^{h(\zeta)}.$$
Since $\beta\geq \alpha$ was arbitrary, we conclude that $\bigcap \limits_{\zeta\in \dom(h)}A_\zeta^{h(\zeta)}$ is bounded  by $\alpha$.
\end{proof} 
\end{theorem}
In the particular case in which $\kappa=\omega$, we have that $\lambda^{<\omega}=\lambda$ for each $\lambda$. Hence, we can rewrite the theorem above as follows:
\begin{theorem} Let $\lambda\leq cof(\mathfrak{c})$. If $\mathfrak{r}(\omega,\lambda)=\mathfrak{c}$, then there is a barely independent family over $\mu$ of size $\lambda.$
\end{theorem}

Assuming that the continuum is a regular cardinal, Theorem \ref{barelyimpliesr} and \ref{rimpliesbarely} yield an equivalence. In particular, since the inequality $\mathfrak{r}<\mathfrak{c}=cof(\mathfrak{c})$ is consistent. The following corollary implies that the existence of barley independent families over $\mathfrak{c}$ is independent from ZFC. 

\begin{corollary}\label{equivalencereapingindependent} Assume that $\mathfrak{c}=cof(\mathfrak{c})$ and let $\kappa\leq \mathfrak{c}$ be another regular cardinal. The followings are equivalent:
\begin{enumerate}
    \item There is a barely independent family over $\mathfrak{c}$ size $\kappa$.
    \item $\mathfrak{r}(\omega\,,\kappa)=\mathfrak{c}.$
\end{enumerate}    
\end{corollary}
If $\kappa \geq \mathfrak{c}$, then any generic extension obtained after adding $\kappa$-many Cohen reals will satisfy that $\mathfrak{r}(\omega,\mathfrak{c})=\mathfrak{c}.$ Thus, the results discussed above apply here. As we see now, there is more to say in this situation. We remark that the theorem below can easily be adapted to the case in which we add Random reals. We leave the details for the reader.

\begin{theorem}\label{theoremcohen}Let $\kappa$ and $\lambda$ be two infinite cardinals and consider $\mathbb{P}=\text{Fin}(\kappa\times \lambda,2)$. Then $$\mathbb{P}\forces{There is a barely independent family over $\kappa$ of size $\lambda$}.$$
\begin{proof}Let $G$ be a $\mathbb{P}$-generic filter and $f_G=\bigcup G$. In $V[G]$, define $A^G_\alpha:=\{\xi\in \kappa \,:\,p(\xi,\alpha)=1\}$ for each $\alpha\in \kappa$. It should be clear that  $\bigcap \mathcal{A}^h\in [\kappa]^\kappa$ for each non-empty finite $\mathcal{A}\subseteq \mathcal{I}^G$.  We will show that $\mathcal{I}^G:=\{A_\alpha^G\,:\,\alpha\in \lambda\}$ is barely independent. For this, let $\mathcal{A}\in [\lambda]^\omega$ and $h:\mathcal{A}\longrightarrow 2$. In $V$, consider $\dot{\mathcal{A}}$ and $\dot{h}$  two $\mathbb{P}$-names for $\mathcal{A}$ and  $h$ respectively, and let $M$ be a countable elementary submodel of a large enough $H(\theta)$ such that $\dot{h}, \dot{\mathcal{A}},\mathbb{P}\in M$. We claim that, in $V[G]$, $\bigcap\limits_{\alpha\in \mathcal{A}}(I^G_\alpha)^{h(\alpha)}\subseteq M$. Indeed, let $p\in \mathbb{P}$ and $\xi\in \kappa\backslash M$. Since $\dot{A}$ is forced to be infinite, the domain of $p$ is finite and both $\dot{A}$ and $\dot{h}$ live in $M$, we can find $\alpha\in (\lambda\cap M)\backslash dom_0(p)$, $i\in 2$ and $q\in \mathbb{P}\cap M$ such that $q\leq p\cap M$ and $q\forces{$\alpha\in \dot{A}$ and $\dot{h}(\alpha)=i$.}$  We now put $r=p\cup q\cup \{((\alpha,\xi), 1-i)\}$. It should be clear that $r$ is well-defined, $r\leq p,q$, and $r\forces{$\xi\not\in (I^G_\alpha)^{\dot{h}(\alpha)}$.}$ This finishes the proof. 
\end{proof}
\end{theorem}
There are two important consequences of the previous theorem.
\begin{corollary}\label{coroimportante}
    
Assume that GCH holds in $V$ and let $\lambda$ be a cardinal of uncountable cofinality.  There is a forcing extension of $W$ of $V$ satisfying following properties:
\begin{itemize}
    \item $2^\kappa=2^\omega=\lambda$ for each uncountable cardinal $\kappa<\lambda$.
    \item For each uncountable cardinal $\kappa\leq \lambda$, there is a barely independent family over $\kappa$ of size $\lambda.$ In particular, every uniform ultrafilter over $\kappa$ is both Tukey top and has character $2^\kappa$.
\end{itemize}    
\end{corollary}
The second consequence is that the existence of a barely independent family over a cardinal  does not yield the existence of a barely independent family over its succesor.
\begin{corollary}\label{lastcoro}It is consistent that $\mathfrak{c}=\omega_2>\omega_1=\mathfrak{r}$ and that there is abarely independent family over $\omega_1$ of size $\omega_1.$
\begin{proof}We just start with a model in which $\mathfrak{r}=\omega_1<cof(\mathfrak{c})\leq\mathfrak{c}$. After adding $\omega_1$-many Cohen reals, the three cardinals mentioned above will have the same values. Hence, there are no Barely independent families over $\mathfrak{c}$ due to Theorem \ref{barelyimpliesr}. On the other hand, there will be a barely independent family over $\omega_1$ of size $\omega_1$ due to Theorem \ref{theoremcohen}.
\end{proof}
\end{corollary}
Our next goal is to show that the existence of barely independent families over $\omega_1$ is also independent of ZFC. Note that the argument provided above can not be used in this case because of the simple reason that $\omega_1\leq \mathfrak{r}(\omega,\omega_1)$. The way of getting this independence is by showing that there are no such families under $\mathfrak{p}=\mathfrak{c}$. This result is unexpected because since $\mathfrak{r}(\omega,\mathfrak{c})=\mathfrak{c}=cof(\mathfrak{c})$ under MA, this axiom does imply the existence of a barely independent family due to Corollary \ref{equivalencereapingindependent}. 

\begin{theorem}\label{teoprinma}Let $\kappa<\mathfrak{p}$ be an uncountable cardinal and let $\mathcal{I}$ be an infinite family of subsets of $\kappa$. Then there is $i\in 2$ and $\mathcal{A}\in [\mathcal{I}]^{\omega}$ for which $|\bigcap\limits_{A\in \mathcal{A}} A^i|=\kappa$.
\begin{proof}Assume towards a contradiction that there is an infinite $\mathcal{I}\subseteq \mathscr{P}(\kappa)$ so that $|\bigcap\limits_{A\in \mathcal{A}}A^i|<\kappa$ for each $i\in 2$ and every  $\mathcal{A}\in [\mathcal{I}]^\omega$. Without loss of generality we may assume that $\mathcal{I}$ is countable. Given $\alpha\in \kappa$ and $i\in 2$, let us define $$\mathcal{O}^i_\alpha:=\{A\in \mathcal{I}\,:\,\alpha\in A^i\}.$$
Observe that $\mathcal{I}$ is equal to the disjoint union of $\mathcal{O}^0_\alpha$ and $\mathcal{O}^1_\alpha.$  Let us fix now a non-principal ultrafilter $\mathcal{U}\subseteq [\mathcal{I}]^\omega$. According to the previous observation, for each $\alpha$ there is $i_\alpha\in 2$ so that $\mathcal{O}^i_\alpha\in \mathcal{U}$. By the pigeonhole principle, we know that there is $X'\in [\kappa]^\kappa$ and  a fixed $i\in 2$ so that $i=i_\alpha$ for any $\alpha\in X'$.  Let us take a countable elementary submodel $M$ of a large enough part of the universe so that $X',\mathcal{I},\kappa\in M$. Now, let $X:=X'\backslash M$. The following fact follows directly from elementarity:\\\\
\noindent
\underline{Fact 1}: Let $\mathcal{B}\in [\mathcal{I}]^{<\omega}$. If there is a non-empty $a\in [X]^{<\omega}$ so that $\mathcal{B}\subseteq \bigcap\limits_{\alpha} \mathcal{O}^i_\alpha$,  then $$\{\min(a)\,:\,a\in [X]^{<\omega}\backslash\{\emptyset\} \text{ and } \mathcal{B}\subseteq \bigcap\limits_{\alpha} \mathcal{O}^i_\alpha\}$$
is unbounded in $\kappa.$\\

We now proceed to define a the $\mathbb{P}$ as the set of all ordered pairs $p=(\mathcal{A}_p,Y_p)$ so that:\begin{itemize}
    \item $Y_p\in [X]^{<\omega}$ and \item $\mathcal{A}_p\subseteq \bigcap\limits_{\alpha\in Y_p}\mathcal{O}^i_\alpha.$ 
    \end{itemize}
    As expected, we put $p\leq q$ whenever $\mathcal{A}_q\subseteq \mathcal{A}_p$ and $Y_q\subseteq Y_p$.\\\\
\noindent
\underline{Claim 1}: $\mathbb{P}$ is $\sigma$-centered.
\begin{claimproof}Given $\mathcal{B}\in [\mathcal{I}]^{<\omega}$, let us define $\mathcal{C}_\mathcal{B}:=\{p\in \mathbb{P}\,:\,\mathcal{A}_p=\mathcal{B}\}$. Since $\mathcal{I}$ is countable, there are countably many $\mathcal{C}_\mathcal{B}$'s. Furthermore, if $p,q\in \mathcal{B}$, it should be clear that $r=(\mathcal{B}, Y_p\cup Y_q)$ is a condition living $\mathcal{C}_\mathcal{B}$ an which is smaller than both $p$ and $q.$ Hence, each  $\mathcal{C}_\mathcal{B}$ is centered.
\end{claimproof}
\noindent
\underline{Claim 2}: Given $\alpha\in \kappa$, let $\mathcal{D}_\alpha=\{p\in \mathbb{P}\,:\,Y_q\backslash \alpha\not=\emptyset\}$. Then $\mathcal{D}_\alpha$ is dense.
\begin{claimproof} Indeed, due to Fact 1, given $p\in \mathbb{P}$, it will happen that $\{\min(Z_q)\,:q\in \mathbb{P}\text{ and }\,\mathcal{A}_q=\mathcal{A}_p\}$ is unbounded in $\kappa$. In particular, this means that there is  $q\in \mathcal{C}_{\mathcal{A}_p}$ with $Z_q\backslash \alpha\not=\emptyset$. We concluded that  $(\mathcal{A}_p, Z_p\cup Z_q)$ is a condition smaller than $p$ which belongs to $\mathcal{D}_\alpha.$
\end{claimproof}
\noindent
\underline{Claim 3}: Given $\mathcal{B}\in [\mathcal{I}]^{<\omega}$, let $\mathcal{E}_{\mathcal{B}}=:\{p\in \mathbb{P}\,:\,\mathcal{A}_p\backslash \mathcal{B}\not=\emptyset\}$. Then $\mathcal{E}_\mathcal{B}$ is dense.
\begin{claimproof}
Indeed, given $p\in \mathbb{P}$ we have that $\bigcap \limits_{\alpha\in Z_p}\mathcal{O}^i_\alpha$ is infinite because each $\mathcal{P}^i_\alpha$ belongs to $\mathcal{U}$. Since $\mathcal{B}$ is infinite, we can take an element $E$ in the previous intersection which do not belong to $\mathcal{B}$. 
The condition $(\mathcal{A}_p\cup \{\mathcal{E}\},Z_p)$ witness density below $p.$
\end{claimproof}
as the family $\{\mathcal{D}_\alpha\,:\alpha\in \kappa\}\cup \{\mathcal{E}_\mathcal{B}\,:\,\mathcal{B}\in [\mathcal{I}]^{<\omega}\}$ has size $\kappa<\mathfrak{p}$ and $\mathbb{P}$ is $\sigma$-centered. There is a filter $G$ intersecting each of the previously mentioned dense sets. In order to finish, just put $\mathcal{A}=\bigcup\limits_{p\in \mathbb{P}} \mathcal{A}_p$ and $Y=\bigcup\limits_{p\in \mathbb{P}}Y_p$. By genericity, $|Y|=\kappa$ and $\mathcal{A}$ is infinite. Furthermore, by definition of $\mathbb{P}$  we have that $\mathcal{A}\subseteq \bigcap\limits_{\alpha\in Y} \mathcal{O}^i_\alpha$. In other words, $$Y\subseteq \bigcap\limits_{A\in \mathcal{A}}A^i.$$
\end{proof}

\end{theorem}
Note that the result above can not be generalized by asking both $\mathcal{I}$ and $\mathcal{A}$ to be of size $\omega_1$. Indeed, the family $\{\omega_1\backslash \alpha\,:\,\alpha\in \omega_1\}$ is a clear counterexample. Additionally, observe that the variable $i$ is also necessary for the theorem. This is because if $\mathcal{I}=\{I_n\,:\,n\in \omega\}$ is $\subseteq$-deacreasing and has empty intersection, any infinite subfamily of $\mathcal{I}$ satisfies the same property.\\
As each barely independent family serves as a counterexample to the conclusion of Theorem \ref{teoprinma}, we conclude the following:
\begin{corollary} Under $\mathfrak{p}>\omega_1$ there are no barely independent families over $\omega_1.$ 
\end{corollary}
Compare the following corollary with Corollary \ref{lastcoro}.
\begin{corollary}Assume that MA holds. Then $\mathfrak{c}$ is the only cardinal for which there are barely independent families.
\end{corollary}
As $\mathfrak{r}(\omega,\mathfrak{c})=\mathfrak{c}$ under CH, we conclude via Theorem \ref{teoprinma} and Theorem \ref{rimpliesbarely} that PID is independent from the existence of barely independent famlies over $\omega_1.$

\section{Open problems}\label{openproblemssection}

We know that the existence of barely independent families over both $\mathfrak{c}$ and $\omega_1$ is independent from ZFC and we know that there can be such families in both cardinals at the same time for arbitrary (consistent) values of $\mathfrak{c}$. By combining the technics for showing the consistency of the non-exitence, we can  have models without barely independent families over both $\mathfrak{c}$ and $\omega_1$. For this,  we just need to look at models in which $\omega_1<\mathfrak{p}\leq \mathfrak{r}<cof(\mathfrak{c}).$ Ths can not be achieve in models in which $\mathfrak{c}=\omega_2.$ Thus, a natural question is:

\begin{problem}Assume that $\mathfrak{c}=\omega_2$. Is it true that there is always a barely independent family?
\end{problem}
A natural way of solving this first problem in the positive would be to show that the following question has also a positive solution.
\begin{problem} Does $\mathfrak{r}=\omega_1$ imply the existence of a barely independent family over $\omega_1?$
\end{problem}
 In general, it would be interesting to analyze the spectrum of existence of $\kappa$-barely independent families for distinct cardinals.\\

In private comunication, Todor\v{c}evi\'c presented me a proof of the fact that the existence of a Luzin set of reals implies the existence of a barely  independent family over $\omega_1$ of size $\omega_1.$ This leaves the question:

\begin{problem}Is it consistent that there is a barely independent family over $\omega_1$ of size $\omega_1$ but there are no Luzin spaces?
\end{problem}

\bibliographystyle{plain}
\bibliography{bibliografia}
{\color{white}hlj}\\\\
\noindent
Jorge Antonio Cruz Chapital\\
Department of Mathematics, University of Toronto, Canada\\
cruz.chapital at utoronto.ca\\
\end{document}